\documentstyle[twoside,amstex,amscd]{article}
\input amssym.def
\input amssym.tex

\def\bc{\begin{center}}
\def\ec{\end{center}}
\def\no{\noindent}
\def\hang{\hangindent\parindent}
\def\textindent#1{\indent\llap{[#1]\enspace}\ignorespaces}
\def\re{\par\hang\textindent}
\begin{document}
\baselineskip12pt \thispagestyle{empty} \vspace*{-3.5mm}
\no\hspace*{88.7mm}\begin{picture}(23.7,13)(0,0)
\setlength{\unitlength}{1mm}
\multiput(0,0)(0,.05){3}{\line(1,0){26.3}}
\multiput(0,13)(0,.05){3}{\line(1,0){26.3}}
\put(0,9.5){\bf\slshape Algebra} \put(0,5.5){\bf\slshape
Colloquium} \put(0,2){\footnotesize \copyright\hspace{.1mm} AMSS
CAS 2004}
\end{picture}

\vspace*{-15mm} \no {\small{\sl Algebra Colloquium} {\bf 11}:4
(2004) 467-476} \vspace*{46mm} \pagestyle{myheadings} \markboth
{\hfill {\small\sl Z.M. Tang}} {{\small\sl Local Homology and
Local Cohomology}\hfill} \vspace*{-1.5 true cm}

\bc{\large\bf Local Homology and Local Cohomology${}^*$}\ec
\bc{{\bf Zhongming Tang}\\
{\small\sl Department of Mathematics, Suzhou University\\
Suzhou 215006, China\\
E-mail: zmtang@@suda.edu.cn}}
\ec \bc{\small Received 5 March 2002\\[1mm]
Communicated by Jin-gen Yang}\ec \vskip0mm

\no\small{{\bf Abstract.}\ \ Let $(R, {\frak m})$ be a local ring,
$I$ a proper ideal of $R$ and $M$ a finitely generated $R$-module
of dimension $d$. We discuss the local homology modules of
$H^d_I(M)$. When $M$ is Cohen-Macaulay, it is proved that
$H^d_{{\frak m}}(M)$ is co-Cohen-Macaulay of N.dimension $d$ and
$H^{\underline{x}} _d(H^d_{{\frak m}}(M))\cong\widehat{M}$ where
$\underline{x}=(x_1,\ldots,x_d)$ is a system of parameters for
$M$.\vskip1mm

\no{\bf 2000 Mathematics Subject Classification:} 13C14,
13D45\vskip1mm

\no{\bf Keywords:}\ \ local homology, local cohomology,
Cohen-Macaulay modules}

\begin{figure}[b]
\vspace{-3mm}
\rule[-2.5truemm]{5cm}{0.1truemm}\\[2mm]
{\small ${}^*$Supported by the National Natural Science Foundation
of China (Grant No. 10071054).}
\end{figure}
\setcounter{page}{467}

\vskip5mm \normalsize

\no{\bf 1\ \ Introduction}\vskip2mm

\no Let $(R,{\frak m})$ be a (commutative Noetherian) local ring
and $M$ a finitely generated $R$-module of dimension $d$. By top
local cohomology modules, we mean the local cohomology modules
$H^d_I(M)$, where $I$ is a proper ideal of $R$. As is well-known,
$H^d_I(M)$ is an Artinian $R$-module (cf. [6]). For
 the Artinian modules over a commutative quasi-local (not necessarily
Noetherian) ring, there is a theory of local homology (cf. [11]).
In this paper, we will mainly discuss the local homology modules
of top local cohomology modules of Cohen-Macaulay modules.

The N.dimension and width of Artinian modules were defined in [8] and [9].
In section 2, we show that N.dim$_R(H^d_I(M))\leq d$, and width$_R(H^d_I(M))
\geq\mbox{min}\{2,d\}$ if $H^d_I(M)\not=0$. For an Artinian module $X$ over
a commutative quasi-local ring $A$, it is always the case that width$_A(X)\leq
\mbox{N.dim}_A(X)$. When the equality holds, we say that $X$ is
co-Cohen-Macaulay. Then it turns out that $H^d_I(M)$ is co-Cohen-Macaulay
if $H^d_I(M)\not=0$ and $d\leq 2$. When $M$ is Cohen-Macaulay, it is proved  that
$M$-regular sequences are $H^d_{{\frak m}}(M)$-coregular. As a corollary, we
have that $H^d_{{\frak m}}(M)$ is co-Cohen-Macaulay of N.dimension $d$ if $M$
is Cohen-Macaulay.  Under the assumption of Cohen-Macaulayness on $M$, the
 main result in section 3 states that $H^{\underline{x}}
_d(H^d_{{\frak m}}(M))\cong\widehat{M}$ where
$\underline{x}=(x_1,\ldots,x_d)$ is a system of parameters for
$M$, which, with another result in this section, demonstrates that
local homology and local cohomology are dual to each other in some
sense. \vskip5mm

\no{\bf 2\ \ Coregular  Sequences on Top Local Cohomology
Modules}\vskip2mm

\no Let $(A,{\frak n})$ be a commutative quasi-local ring and $X$
an Artinian $A$-module. Roberts [9] introduced a dimension for
$X$. Following Kirby [3], we call this dimension as N.dimension
(Noetherian dimension) of $X$, denoted by N.dim$_A(X)$.
N.dim$_A(X)$ is defined as follows: we put $\mbox{N.dim}_A(X)=-1$
when $X=0$, then, inductively, let $r\geq 0$ be an integer, when
$\mbox{N.dim}_A(X)<r$ is false and for any ascending chain
$$
X_0\subseteq X_1\subseteq X_2\subseteq\cdots
$$
of submodules of  $X$ there exists an integer $n$ such that
$\mbox{N.dim}_A(X_ {m+1}/X_m)<r$ for all  $m\geq n$, we put
N.dim$_A(X)=r$. Then, N.dim$_A(X)=0$ if and only if $X$ has finite
length and N.dim$_A(X)$ is the least integer $r$ such that
$0:_X(x_1,\ldots,x_r)$ has finite length for some
$x_1,\ldots,x_r\in{\frak n}$ (cf. [9]).

Let $x_1,\ldots,x_n\in{\frak n}$. We say that $x_1,\ldots,x_n$ is an
$X$-coregular sequence if
$$
0:_X(x_1,\ldots,x_{i-1})\stackrel{x_i}{\longrightarrow}
0:_X(x_1,\ldots, x_{i-1})
$$
is surjective for $i=1,\ldots,n$. The width of $X$, denoted by width$_A(X)$,
 is the length of a maximal $X$-coregular sequence in ${\frak n}$. For any
$X$-coregular element $x\in{\frak n}$, we have that
$$
\mbox{N.dim}_A(0:_Xx)=\mbox{N.dim}_A(X)-1
$$
and
$$
\mbox{width}_A(0:_Xx)=\mbox{width}_A(X)-1,
$$
(cf. [8] and [9]).

Let $M$ be a Cohen-Macaulay module  of dimension $d$ over a local
ring, the following proposition states that  regular sequences on
$M$ are coregular on the top local cohomology module $H_{{\frak
m}}^d(M)$.\vskip3mm

\no{\bf Proposition 2.1.}\ \ {\it Let $(R,{\frak m})$ be a local
ring and $M$ a Cohen-Macaulay $R$-module. If $x_1,\ldots,x_i\in
{\frak m}$ is a regular sequence on $M$, then it is $H_{{\frak
m}}^d(M)$-coregular. Furthermore, for any $n\geq 1$, there is an
isomorphism
$$
\alpha_{i,n}: H^{d-i}_{{\frak
m}}(M/(x_1^n,\ldots,x_i^n)M)\longrightarrow 0:_{H_{{\frak
m}}^d(M)} (x_1^n,\ldots,x_i^n),
$$
and, for all $n\geq 1$, the following diagram is commutative
$$\CD
H^{d-i}_{{\frak m}}(M/(x_1^{n+1},\ldots,x_i^{n+1})M)
@>\alpha_{i,n+1} >>
0:_{H^d_{{\frak m}}(M)}(x_1^{n+1},\ldots,x_i^{n+1})\\
@V \nu_{n} VV @VV x_1\cdots x_i V \\
H^{d-i}_{{\frak m}}(M/(x_1^{n},\ldots,x_i^{n})M)@>
\alpha_{i,n}>>0:_{H^d_{{\frak m}}(M)}(x_1^{n},\ldots,x_i^{n})
\endCD$$
where $\nu_n$ is the natural homomorphism.}\vskip2mm

\no{\it Proof.}\ \ We use induction on $i$. Set
$\underline{x}^n_i=(x_1^n,\ldots,x_i^n)$. Considering the long
exact sequences of local cohomology modules of the diagram
$$\CD
0@>>> M @> x_1^{n+1} >> M @>>> M/x_1^{n+1}M @>>>0\\
&& @Vx_1VV @VVV @VVV\\
0@>>> M @>x_1^n >> M @>>> M/x_1^nM @>>> 0
\endCD$$
where the middle vertical homomorphism is the identity map and the
right vertical homomorphism is the natural map, as $M/x_1^nM$ is
Cohen-Macaulay of dimension $d-1$, we get the following
commutative diagram with exact rows
$$\CD
0@>>>H_{{\frak m}}^{d-1}(M/x_1^{n+1}M)@>>> H^d_{{\frak
m}}(M)@>x_1^{n+1}>>H^d_{{\frak m}}(M)@>>>0\\
&& @VVV @Vx_1VV @VVV &&\\
0@>>>H^{d-1}_{{\frak m}}(M/x_1^nM)@>>>H^d_{{\frak m}}(M)
@>x^n_1>>H^d_{{\frak m}}(M)@>>> 0
\endCD$$
Thus  $x_1$ is $H_{{\frak m}}^d(M)$-coregular and, for any $n\geq 1$, there is an isomorphism
$$
\alpha_{1,n}: H^{d-1}_{{\frak m}}(M/x_1^nM)\longrightarrow
0:_{H^d_{{\frak m}}(M)}x_1^n,
$$
and, for all $n\geq 1$, the following diagram is commutative
$$\CD
H^{d-1}_{{\frak m}}(M/x_1^{n+1}M)@>\alpha_{1,n+1}>>
0:_{H^d_{{\frak m}}(M)}x_1^{n+1}\\
@V\nu_{1,n}VV @VV x_1 V\\
H^{d-1}_{{\frak m}}(M/x_1^nM)@>\alpha_{1,n}>> 0:_{H^d_{{\frak
m}}(M)}x_1^n
\endCD
$$
so, the result is true for $i=1$.

Now suppose that $i>1$ and the result is true for $i-1$. Then, we have the
following commutative diagram
$$\CD
0:_{H_{{\frak m}}^{d-i+1}(M/\underline{x}_{i-1}^{n+1}M)}x_i^{n+1}
@>\overline{\alpha}_{i-1,n+1}>>
0:_{H^d_{{\frak m}}(M)}\underline{x}_i^{n+1}\\
@V \overline{\nu}_{i-1,n}VV @VV x_1\cdots x_{i-1}V\\
0:_{H_{{\frak m}}^{d-i+1}(M/\underline{x}_{i-1}^nM)}x_i^{n+1}
@>\alpha'_{i-1,n}>>0:_{H^d_{{\frak m}}(M)}(\underline{x}_{i-1}^n,x_i^{n+1})\\
@V x_i VV @VV x_iV\\
0:_{H_{{\frak m}}^{d-i+1}(M/\underline{x}_{i-1}^nM)}x_i^n
@>\overline{\alpha}_{i-1,n}>> 0:_{H^d_{{\frak
m}}(M)}\underline{x}_i^n
\endCD
$$
where $\overline{\alpha}_{i-1,n+1},\alpha'_{i-1,n}$ and
$\overline{\alpha}_{i-1,n}$ are the isomorphisms obtained by restricting
$\alpha_{i-1,n+1},\alpha_{i-1,n}$ and $\alpha_{i-1,n}$ on the corresponding
submodules respectively, and $\overline{\nu}_{i-1,n}$ is also the restriction
of $\nu_{i-1,n}$.
 Considering the long exact sequences of local cohomology modules
of the following diagram
$$\CD
0@>>> M/\underline{x}_{i-1}^{n+1}M @>x_i^{n+1}>>
M/\underline{x}_{i-1}^{n+1}M@>>> M/\underline{x}_i^{n+1}M@>>> 0\\
&& @VVV @VVV @VVV\\
0@>>> M/\underline{x}_{i-1}^nM@>x_i^{n+1}>>
M/\underline{x}_{i-1}^nM@>>>
M/\underline{x}_{i-1,i}^{n,n+1}M@>>> 0\\
&& @V x_i VV @VVV @VVV\\
0@>>> M/\underline{x}_{i-1}^nM@>x_i^n>>
M/\underline{x}_{i-1}^nM@>>> M/\underline{x}_i^nM@>>> 0
\endCD
$$
where
$\underline{x}_{i-1,i}^{n,n+1}=(\underline{x}_{i-1}^n,x_i^{n+1})$
and unlabeled maps are natural maps, we have the following
commutative diagram with exact rows by using temporarily $H^i_N$
to denote $H^i_{{\frak m}}(N)$
$$\CD
0@>>>H^{d-i}_{M/\underline{x}_i^{n+1}M}@>>>
H^{d-i+1}_{M/\underline{x}_{i-1}^{n+1}M} @>x_i^{n+1}>>
H^{d-i+1}_{M/\underline{x}_{i-1}^{n+1}M}@>>>0\\
&& @VVV @VVV @VVV\\
0@>>>H^{d-i}_{M/\underline{x}_{i-1,i}^{n,n+1}M}@>>>
H^{d-i+1}_{M/\underline{x}_{i-1}^nM} @>x_i^{n+1}>>
H^{d-i+1}_{M/\underline{x}_{i-1}^nM}@>>>0\\
&&@VVV @Vx_iVV @VVV\\
0@>>> H^{d-i}_{M/\underline{x}_i^nM}@>>>
H^{d-i+1}_{M/\underline{x}_{i-1}^nM}
@>x_i^n>>H^{d-i+1}_{M/\underline{x}_{i-1}^nM}@>>> 0.
\endCD
$$
Then $x_iH_{{\frak m}}^{d-i+1}(M/\underline{x}_{i-1}M)=H_{{\frak
m}}^{d-i+1}(M/\underline{x}_{i-1}M)$. As $\alpha_{i-1,1}$ is an
isomorphism and $x_1,\ldots,x_{i-1}$ is $H^d_{{\frak
m}}(M)$-coregular, we see that $x_1,\ldots,x_i$ is an $H^d_{{\frak
m}}(M)$-coregular sequence. Further,  we get the following
commutative diagram
$$\CD
H_{{\frak m}}^{d-i}(M/\underline{x}_i^{n+1}M)@>>>
0:_{H_{{\frak m}}^{d-i+1}(M/\underline{x}_{i-1}^{n+1}M)}x_i^{n+1}\\
@VVV @VV\overline{\nu}_{i-1,n}V\\
H_{{\frak m}}^{d-i}(M/(\underline{x}_{i-1}^n,x_i^{n+1})M)@>>>
0:_{H_{{\frak m}}^{d-i+1}(M/\underline{x}_{i-1}^nM)}x_i^{n+1}\\
@VVV @VVx_iV\\
H_{{\frak m}}^{d-i}(M/\underline{x}_i^nM)@>>>
0:_{H_{{\frak
m}}^{d-i+1}(M/\underline{x}_{i-1}^nM)}x_i^n
\endCD
$$
where the rows are isomorphisms and the composition of the two
left vertical maps is just the homomorphism reduced from the
natural homomorphism $M/\underline{x}_i^{n+1}M \longrightarrow
M/\underline{x}_i^nM$. Thus, for any $n\geq 1$, we have an
isomorphism
$$
\alpha_{i,n}:H_{{\frak
m}}^{d-i}(M/\underline{x}_i^nM)\longrightarrow 0:_{H^d_{{\frak
m}}(M)}\underline{x}_i^n,
$$
and, for all $n\geq 1$ the following diagram is commutative
$$\CD
H_{{\frak m}}^{d-i}(M/\underline{x}_i^{n+1}M)@>\alpha_{i,n+1}>>
0:_{H^d_{{\frak m}}(M)}\underline{x}_i^{n+1}\\
@V\nu_{i,n}VV @VV x_1\cdots x_iV\\
H_{{\frak m}}^{d-i}(M/\underline{x}_i^nM)@>\alpha_{i,n}>>
0:_{H^d_{{\frak m}}(M)}\underline{x}_i^n.
\endCD
$$
The proposition follows. \hfill$\Box$\vskip2mm

The following corollary will be used in the next section.\vskip3mm

\no{\bf Corollary 2.2.}\ \ {\it Let $(R,{\frak m})$ be a local
ring, $M$ a finitely generated $R$-module of dimension $d$ and
$x_1,\ldots,x_d$ a system of parameters for $M$. Suppose that $M$
is Cohen-Macaulay. Then, for any $n\geq 1$, there is an
isomorphism
$$
\alpha_n: M/(x_1^n,\ldots,x_d^n)M\longrightarrow 0:_{H_{{\frak
m}}^d(M)} (x_1^n,\ldots,x_d^n),
$$
and, for all $n\geq 1$, the following diagram is commutative
$$\CD
M/(x_1^{n+1},\ldots,x_d^{n+1})M@>\alpha_{n+1}>>
0:_{H^d_{{\frak m}}(M)}(x_1^{n+1},\ldots,x_d^{n+1})\\
@V \nu_nVV @VV x_1\cdots x_d V\\
M/(x_1^n,\ldots,x_d^n)M @>\alpha_n>> 0:_{H^d_{{\frak
m}}(M)}(x_1^n,\ldots,x_d^n)
\endCD
$$
where $\nu_n$ is the natural homomorphism.} \vskip3mm

For the N.dimension and width of top local cohomology modules, we
have the following two propositions.\vskip3mm

\no{\bf Proposition 2.3.}\ \ {\it  Let $(R,{\frak m})$ be a local
ring, $M\not=0$ a finitely generated $R$-module of dimension $d$
and $I$ a proper ideal of $R$. Then
$$
\mbox{N.dim}_R(H_I^d(M))\leq d.
$$}

\no{\it Proof.}\ \  We use induction on $d$. When $d=0$, $M$ has
finite length. Then $M$ is annihilated by some power of $I$, hence
$H_I^0(M)=M$ . Thus $\mbox{N.dim}_R(H_I^0(M))\\ =0$.

Now assume that $d>0$.
As $H_{{\frak m}}^0(M)$ is a submodule of $M$ which has finite length,
we see that $M/H_{{\frak m}}^0(M)\not=0$ is also of dimension $d$. Further,
as $H^i_I(H_{{\frak m}}^0(M))=0$ for all $i>0$, from the long exact sequence of local
cohomology modules of the short exact sequence
$$
0\rightarrow H_{{\frak m}}^0(M)\rightarrow M\rightarrow M/H_{{\frak m}}^0(M)\rightarrow 0
$$
we get that
$$
H^d_I(M)\cong H^d_I(M/H_{{\frak m}}^0(M)).
$$
Thus, we may assume that $H_{{\frak m}}^0(M)=0$. Then, in ${\frak
m}$, there exists non-zero divisors on $M$. Let $x\in{\frak m}$ be
any non-zero divisor on $M$. Thus dim$(M/xM)=d-1$, hence
$H_I^d(M/xM)=0$. Considering the long exact sequence of local
cohomology modules of the short exact sequence
$$
0\rightarrow M\stackrel{x}{\rightarrow}M\rightarrow M/xM\rightarrow 0,
$$
we have the following exact sequence
$$
\cdots\cdots\rightarrow H_I^{d-1}(M/xM)\rightarrow H_I^d(M)\stackrel{x}
{\rightarrow}H_I^d(M)\rightarrow 0.
$$
Then $xH_I^d(M)=H_I^d(M)$ and $0:_{H^d_I(M)}x$ is isomorphic to a
factor module of $H^{d-1}_I(M/xM)$, hence
$$
\mbox{N.dim}_R(H^d_I(M))=\mbox{N.dim}_R(0:_{H^d_I(M)}x)+1
$$
and
$$
\mbox{N.dim}_R(0:_{H^d_I(M)}x)\leq\mbox{N.dim}_R(H^{d-1}_I(M/xM)).
$$
Therefore, by induction assumption, we have
\begin{eqnarray*}
&&
\mbox{N.dim}_R(H^d_I(M))\\
&=&
\mbox{N.dim}_R(0:_{H^d_I(M)}x)+1\\
&\leq&
\mbox{N.dim}_R(H^{d-1}_I(M/xM))+1\\
&\leq&
(d-1)+1\\
&=&
d
\end{eqnarray*}
The proposition follows.\hfill$\Box$\vskip2mm

\no{\bf Proposition 2.4.}\ \ {\it Let $(R,{\frak m})$  be a local
ring, $M \not=0$ a finitely generated $R$-module of dimension $d$
and $I$ a proper ideal of $R$. If $H_I^d(M)\not=0$, then
$$
\mbox{width}_R(H_I^d(M))\geq\min\{2,d\}.
$$}

\no{\it Proof.} \ \ When $d=0$, it is trivial. We assume that
$d\geq 1$. As in the proof of proposition 2.3, we may assume that
there exists $x\in{\frak m}$ which is a non-zero divisor on $M$,
and we have $xH_I^d(M)=H_I^d(M)$ and an epimorphism
$H_I^{d-1}(M/xM)\longrightarrow 0:_{H_I^d(M)}x$. Then
width$_R(H_I^d(M)) \geq 1$ and as $H_I^d(M)$ is Artinian, we have
$0:_{H_I^d(M)}x\not =0$, hence $H_I^{d-1}(M/xM)\not=0$.

Suppose that $d\geq 2$. Then, by using the same arguments to
$H_I^{d-1}(M/xM)$, we get some $y\in{\frak m}$ such that
$yH_I^{d-1}(M/xM)= H_I^{d-1}(M/xM)$. But, $0:_{H_I^d(M)}x$ is a
homomorphic image of $H_I^{d-1} (M/xM)$, so,
$y(0:_{H_I^d(M)}x)=0:_{H_I^d(M)}x$. Thus, $x,y$ is an
$H_I^d(M)$-coregular sequence, hence, width$_R(H_I^d(M))\geq 2$.
It follows that width$_R(H_I^d(M))\geq \min\{2,d\}$.
\hfill$\Box$\vskip2mm

For an Artinian module $X$ over a commutative quasi-local ring
$A$, it is always the case that width$_A(X)\leq\mbox{N.dim}_A(X)$.
When the equality holds, we say that $X$ is co-Cohen-Macaulay.
Co-Cohen-Macaulay modules were introduced and studied in [11] and
[12].

Let $k$ be a field and $R=k[[X_1,\ldots,X_n]]$ be a formal power
series ring. Then $R$ is a local ring and $(X_1,\ldots,X_n)$ is
the maximal ideal. Let $K=k[X_1^{-1},\ldots,X_n^{-1}]$ be the
$k$-module of polynomials in $X_1^{-1}, \ldots,X_n^{-1}$ with
coefficients in $k$. Then $K$ is an Artinian $R$-module ( see
[2]). It is easy to see that $X_1,\ldots,X_n$ is a $K$-coregular
sequence and $0:_K(X_1,\ldots,X_n)=k$ has finite length, hence,
width$_R(K)=\mbox{N.dim} _R(K)=n$. Thus $K$ is a co-Cohen-Macaulay
$R$-module.

From propositions 2.3 and 2.4, we have immediately the
following\vskip3mm

\no{\bf Corollary 2.5.}\ \ {\it  Let $(R,{\frak m})$ be a local
ring, $M$ a finitely generated $R$-module of dimension $d$ and $I$
a proper ideal of $R$. If $H_I^d(M)\not=0$ and $d\leq 2$, then
$H_I^d(M)$ is co-Cohen-Macaulay.}\vskip2mm

Let $(R,{\frak m})$ be a local ring and $M\not=0$ a Cohen-Macaulay
$R$-module of dimension $d$. Let $x_1,\ldots,x_d$ be a system of parameters for $M$.
 Then, by  proposition 2.1, $x_1,\ldots,x_d$ is $H_{{\frak m}}^d(M)$-coregular, hence,
width$_R(M)\geq d$. Thus, the following proposition follows from
proposition 2.3.\vskip3mm

\no{\bf Proposition 2.6.}\ \ {\it Let $(R,{\frak m})$ be a local
ring and $M\not=0$ a finitely generated $R$-module of dimension
$d$. If $M$ is Cohen-Macaulay, then $H_{{\frak m}}^d(M)$ is
co-Cohen-Macaulay of N.dimension $d$.}\vskip5mm

\no{\bf 3\ \ Local Homology of Top Local Cohomology
Modules}\vskip2mm

\no Let $(A,{\frak n})$ be a commutative quasi-local ring and $X$
an Artinian $A$-module. For any sequence
$\underline{x}=(x_1,\ldots,x_r)$ in ${\frak n}$ and an integer
$i$, the $i$-th local homology module $H_i^{\underline{x}}(X)$ of
$X$ with respect to $\underline{x}$ is defined by
$$
\lim_{\stackrel{\longleftarrow}{n}}\{H_i(K.(x_1^n,\ldots,x_r^n;X))\},
$$
where $K.(x_1^n,\ldots,x_r^n;X)$ is the Koszul complex of $X$ with
respect to $x_1^n,\ldots, \\x_r^n$. Then $H_i^{\underline{x}}(~)$
is an additive, $A$-linear covariant functor from the category of
Artinian $A$-modules and $A$-homomorphisms to the category of
$A$-modules and $A$-homomorphisms, and, for any short exact
sequence of Artinian $A$-modules, there is a long exact sequence
of local homology modules. Furthermore, it is proved in [11] that
$H_i^{\underline{x}}(X)=0$ for $i>\mbox{N.dim}_A(X)$ and $\inf\{i|
H_i^{\underline{x}}(X)\not =0\}$ is the width of $X$ on the ideal
$(x_1,\ldots,x_r)$, i.e., the length of a maximal $X$-coregular
sequence in the ideal $(x_1,\ldots,x_r)$.\vskip2mm

Now we can prove the main result.\vskip3mm

\no{\bf Theorem 3.1.}\ \ {\it Let $(R,{\frak m})$ be a local ring,
$M$ a finitely generated $R$-module of dimension $d$ and
$\underline{x}=(x_1,\ldots,x_d)$  a system of
 parameters for $M$. If $M$ is Cohen-Macaulay, then
 $$
H_i^{\underline{x}}(H^j_{{\frak m}}(M))\cong\left\{
\begin{array}{ll}
\widehat{M}&i=j=d\\
0&\mbox{otherwise}
\end{array}
\right.
$$
where $\widehat{M}$ is the ${\frak m}$-adic completion of
$M$.}\vskip2mm

\no{\it Proof.}\ \ Since $M$ is Cohen-Macaulay, we have that
$H^j_{{\frak m}}(M)=0$ if $j\not=d$. On the other hand, as
$H^d_{{\frak m}}(M)$ is co-Cohen-Macaulay of N.dimension $d$ and
$x_1,\ldots,x_d$ is an $H^d_{{\frak m}}(M)$-coregular sequence by
propositions 2.6 and 2.1, we have that
$H_i^{\underline{x}}(H^d_{{\frak m}}(M))=0$ if $i\not=d$. Hence,
it is enough to show that $H_d^{\underline{x}}(H^d_{{\frak
m}}(M))\cong \widehat{M}$.

Since
$$
H_d^{\underline{x}}(H^d_{{\frak m}}(M))\cong\lim
_{\stackrel{\longleftarrow}{n}}\{0:_{H^d_{{\frak
m}}(M)}(x_1^n,\ldots,x_d^n); \varphi_n\},
$$
where
\begin{eqnarray*}
\varphi_n: 0:_{H^d_{{\frak m}}(M)}(x_1^n,\ldots,x_d^n)
&\longrightarrow& 0:_{H^d_{{\frak
m}}(M)}(x_1^{n-1},\ldots,x_d^{n-1})\\
m&\mapsto& x_1\cdots x_dm
\end{eqnarray*}
and
$$
\widehat{M}\cong\lim_{\stackrel{\longleftarrow}{n}}\{
M/(x_1^n,\ldots,x_d^n)M;\psi_n\},
$$
where $\psi_n:M/(x_1^n,\ldots,x_d^n)M\longrightarrow M/(x_1^{n-1},
\ldots,x_d^{n-1})M$ is the natural homomorphism, and, by corollary
2.2, the two inverse systems are isomorphic, it follows that
$H_d^{\underline{x}}(H^d_{{\frak m}}(M))\cong\widehat{M}$.
\hfill$\Box$\vskip2mm

Because local homology functors and local cohomology functors are
$R$-linear, the following corollary follows from theorem
3.1.\vskip2mm

\no{\bf Corollary 3.2.}\ \ {\it Let $(R,{\frak m})$ be a local
ring and $M$ a finitely generated $R$-module of dimension $d$. If
$M$ is Cohen-Macaulay, then
$$
\mbox{ann}_R(H^d_{{\frak m}}(M))=\mbox{ann}_R(M)\bigskip.
$$}

By using the dual arguments to proposition 2.1, we can show the
following\vskip2mm

\no{\bf Lemma 3.3.}\ \ {\it Let $(A,{\frak n})$ be a commutative
quasi-local ring, $X$ an Artinian $A$-module of N.dimension $d$
and $\underline{x}=(x_1,\ldots,x_d)\subseteq{\frak n}$ such that
$0:_X\underline{x}$ has finite length. Suppose that $X$ is
co-Cohen-Macaulay. Then, for any $n\geq 1$, there is an
isomorphism
$$
\beta_n:H_d^{\underline{x}}(X)/(x_1^n,\ldots,x_d^n)H_d^{\underline{x}}(X)\longrightarrow
0:_X(x_1^n,\ldots,x_d^n),
$$
and, for all $n\geq 1$, the following diagram is commutative
$$\CD
H_d^{\underline{x}}(X)/(x_1^n,\ldots,x_d^n)H_d^{\underline{x}}(X)
@>\beta_n>>0:_X(x_1^n,\ldots,x_d^n)\\
@V x_1\cdots x_d VV @VV i_nV\\
H_d^{\underline{x}}(X)/(x_1^{n+1},\ldots,x_d^{n+1})H_d^{\underline{x}}(X)
@>\beta_{n+1}>> 0:_X(x_1^{n+1},\ldots,x_d^{n+1})
\endCD
$$
where $i_n$ is the inclusion map.}\vskip2mm

The following theorem is dual to theorem 3.1.\vskip3mm

\no{\bf Theorem 3.4.}\ \ {\it Let $(R,{\frak m})$ be a local ring,
$X$ an Artinian $R$-module of N.dimension $d$ and
$\underline{x}=(x_1,\ldots,x_d)\subseteq{\frak m}$ such that
$0:_X\underline{x}$ has finite length. If $X$ is
co-Cohen-Macaulay, then
$$
H^d_{\underline{x}}(H_d^{\underline{x}}(X))\cong X,
$$
where $\underline{x}$ is also interpreted as the ideal generated
by $x_1,\ldots,x_d$.}

\no{\it Proof.}\ \ Since
$$
H^d_{\underline{x}}(H_d^{\underline{x}}(X))\cong\lim_{\stackrel{\longrightarrow}{n}}
\{H_d^{\underline{x}}(X)/(x_1^n,\ldots,x_d^n)H_d^{\underline{x}}(X);\varphi_n\}
$$
where
\begin{eqnarray*}
\varphi_n:H_d^{\underline{x}}(X)/(x_1^n,\ldots,x_d^n)H_d^{\underline{x}}(X)
&\rightarrow&
H_d^{\underline{x}}(X)/(x_1^{n+1},\ldots,x_d^{n+1})H_d^{\underline{x}}(X)\\
m+(x_1^n,\ldots,x_d^n)H_d^{\underline{x}}(X)&\mapsto& x_1\cdots
x_dm+(x_1^{n+1},\ldots,x_d^{n+1})H_d^{\underline{x}}(X),
\end{eqnarray*}
and
$$
X=\bigcup_{n\geq
1}0:_X(x_1^n,\ldots,x_d^n)=\lim_{\stackrel{\longrightarrow}{n}}
\{0:_X(x_1^n,\ldots,x_d^n);i_n\}
$$
where $i_n:0:_X(x_1^n,\ldots,x_d^n)\longrightarrow
0:_X(x_1^{n+1},\ldots,x_d^{n+1})$ is the inclusion map, but, by
lemma 3.3, the two direct systems are isomorphic, it follows that
$H^d_{\underline{x}}(H^{\underline{x}}_d(X))\cong X$.
\hfill$\Box$\vskip2mm

\no{\it Remark 3.5.}\ \ Under the assumption of theorem 3.4, we
see from the proof of lemma 3.3 that $x_1,\ldots,x_d$ is an
$H_d^{\underline{x}}(X)$-regular sequence and
$$
H_d^{\underline{x}}(X)/(x_1,\ldots,x_d)H_d^{\underline{x}}(X)\not=0
$$
has finite length. Hence $H_d^{\underline{x}}(X)$ would be
Cohen-Macaulay if one can show that $H_d^{\underline{x}}(X)$ is
finitely generated. On the other hand, if dim$(R)=d$, then
$H_d^{\underline{x}}(X)$ is a maximal Cohen-Macaulay
$R$-module.\vskip4mm

\parindent=8mm
\no{\bf References} \vskip2mm {\small

\re{1} M.P. Brodmann, R.Y. Sharp, {\it Local Cohomology},
Cambridge Univ. Press, 1998.

\re{2} D. Kirby, Artinian modules and Hilbert polynomials, {\it
Quart. J. Math. (Oxford)} (2) {\bf 24} (1973)  47--57.

\re{3} D. Kirby, Dimension and length for Artinian modules, {\it
Quart. J. Math. (Oxford)} (2) {\bf 41} (1990) 419--429.

\re{4} I.G. MacDonald, Secondary representations of modules over a
commutative ring, {\it Symposia Math.} {\bf 11} (1973) 23--43.

\re{5} I.G. MacDonald,  R. Y. Sharp,  An elementary proof of the
non vanishing of certain local cohomology modules, {\it Qurat. J.
Math. (Oxford) } {\bf 23} (1972) 197--204.

\re{6} L. Melkersson, Some applications of a criterion for
artinianness of a module, {\it J. Pure and Appl. Alg. } {\bf 101}
(1995) 291--303.

\re{7} L. Melkersson, P. Schenzel, The co-localization of an
Artinian module,
 {\it Proc. Edin. Math. Soc.} {\bf 38} (1995) 121--131.

\re{8} A. Ooishi, Matlis duality and width of a module, {\it
Hiroshima Math. J.} {\bf  6} (1976) 573--587.

\re{9} R.N. Roberts, Krull dimension for Artinian modules over
quasi local commutative rings, {\it Quart. J. Math. (Oxford)(3)}
{\bf 26} (1975) 269--273.

\re{10} R.Y. Sharp, Local cohomology theory in commutative
algebra, {\it Quart. J. Math. (Oxford)}(2) {\bf 21} (1970)
425--434.

\re{11} Z.M. Tang, Local homology theory for Artinian modules,
{\it Comm. Alg. } {\bf 22} (1994) 1675--1684.

\re{12} Z.M. Tang,  H. Zakeri, Co-Cohen-Macaulay modules and
modules of generalized fractions, {\it Comm. Alg.} {\bf 22} (1994)
2173--2204.
\end{document}